# Stability Results for the Continuity Equation

**Iasson Karafyllis[*] and Miroslav Krstic[**]**

[*]Dept. of Mathematics, National Technical University of Athens, Zografou Campus, 15780, Athens, Greece,
email: iasonkar@central.ntua.gr; iasonkaraf@gmail.com

[**]Dept. of Mechanical and Aerospace Eng., University of California, San Diego, La Jolla, CA 92093-0411, U.S.A., email: krstic@ucsd.edu

**Abstract**

We provide a thorough study of stability of the 1-D continuity equation, which models many physical conservation laws. In our system-theoretic perspective, the velocity is considered to be an input. An additional input appears in the boundary condition (boundary disturbance). Stability estimates are provided in all $L^p$ state norms with $p>1$, including the case $p=+\infty$. However, in our Input-to-State Stability estimates, the gain and overshoot coefficients depend on the velocity. Moreover, the logarithmic norm of the state appears instead of the usual norm. The obtained results can be used in the stability analysis of larger models that contain the continuity equation. In particular, it is shown that the obtained results can be used in a straightforward way for the stability analysis of non-local, nonlinear manufacturing models under feedback control.

## 1. Introduction

The continuity equation is the conservation law of every quantity that is transferred only by means of convection. It arises in many models in mathematical physics and for the 1-D case takes the form

$$\frac{\partial \rho}{\partial t}+\frac{\partial}{\partial x}(\rho v)=0 \qquad (1.1)$$

where $t$ denotes time, $x$ is the spatial variable, $\rho$ is the density of the conserved quantity and $v$ is the velocity of the medium. Equation (1.1) is used in fluid mechanics (conservation of mass; see Chapter 13 in [26] and [5]), in electromagnetism (conservation of charge; see Chapter 13 in [26]), in traffic flow models (conservation of vehicles; see Chapter 2 in [15], [21,28] and references therein) as well as in many other cases where $\rho$ is not necessarily the density of a conserved quantity (e.g., in shallow water equations the continuity equation is obtained as a consequence of the conservation of mass with the fluid height in place of $\rho$; see [27]). In many cases, the conserved quantity can only have positive values (e.g., mass density, vehicle density) and equation (1.1) comes together with the positivity requirement $\rho>0$.

Although the continuity equation is used extensively in many mathematical models, its stability properties have, surprisingly, not been studied in detail (but see [1,6] for other aspects of the continuity equation). This is largely because the continuity equation usually appears as a part of a larger mathematical model, which also describes the evolution of the velocity profile. In other words, the continuity equation does not have the velocity as an independent input but is accompanied by at least one more equation: a differential equation (momentum balance) in fluid mechanics and traffic flow (see also [13] for oil drilling), or a non-local equation in manufacturing models (see [7,8,9,23]).

In this paper we study the stability properties of the continuity equation on its own. We adopt a system-theoretic perspective, where $v$ is an input. When the velocity profile is given, the continuity equation falls within the framework of transport PDEs, which are studied heavily in the literature (see for instance [2,3,4,10,11,12,14,17,18,19,21,22,24,25,29]). In this framework, the continuity equation is a bilinear transport PDE. Using the characteristic curves and a Lyapunov analysis we are in a position to establish stability estimates that look like Input-to-State Stability (ISS) estimates with respect to all inputs: the input $v$ as well as boundary inputs (Theorem 2.1). The stability estimates are provided in all $L^p$ state norms with $p > 1$, including the case $p = +\infty$. However, the obtained estimates are not precisely ISS estimates since both the gain and overshoot coefficients depend on the input $v$. Moreover, the logarithmic norm of the state appears instead of the usual norm; this is common in many systems where the state variable is positive (see [18]).

The study of the stability properties of the continuity equation can be useful in the stability analysis of larger models (by using small-gain arguments; see [20]). Indeed, we show that the obtained results can be used in a straightforward way for the stability analysis of non-local, nonlinear manufacturing models (Theorem 3.1).

The structure of the paper is as follows. In Section 2, we present the stability estimates for the continuity equation. Section 3 is devoted to the application of the obtained results to non-local, nonlinear manufacturing models. The proofs of all results are provided in Section 4. Finally, the concluding remarks of the present work are given in Section 5.

**Notation:** Throughout the paper, we adopt the following notation.

* $\Re_+ := [0, +\infty)$. Let $u : \Re_+ \times [0,1] \to \Re$ be given. We use the notation $u[t]$ to denote the profile at certain $t \geq 0$, i.e., $(u[t])(x) = u(t,x)$ for all $x \in [0,1]$. $L^p(0,1)$ with $p \geq 1$ denotes the equivalence class of measurable functions $f : [0,1] \to \Re$ for which $\|f\|_p = \left( \int_0^1 |f(x)|^p \, dx \right)^{1/p} < +\infty$. $L^\infty(0,1)$ denotes the equivalence class of measurable functions $f : [0,1] \to \Re$ for which $\|f\|_\infty = \underset{x \in (0,1)}{\operatorname{ess\,sup}} (|f(x)|) < +\infty$. We use the notation $f'(x)$ for the derivative at $x \in [0,1]$ of a differentiable function $f : [0,1] \to \Re$.

* Let $S \subseteq \Re^n$ be an open set and let $A \subseteq \Re^n$ be a set that satisfies $S \subseteq A \subseteq cl(S)$. By $C^0(A; \Omega)$, we denote the class of continuous functions on $A$, which take values in $\Omega \subseteq \Re^m$. By $C^k(A; \Omega)$, where $k \geq 1$ is an integer, we denote the class of functions on $A \subseteq \Re^n$, which takes values in $\Omega \subseteq \Re^m$ and has continuous derivatives of order $k$. In other words, the functions of class $C^k(A; \Omega)$ are the functions which have continuous derivatives of order $k$ in $S = \operatorname{int}(A)$ that can be continued continuously to all points in $\partial S \cap A$. When $\Omega = \Re$ then we write $C^0(A)$ or $C^k(A)$.

* A left-continuous function $f : [0,1] \to \Re$ (i.e. a function with $\lim_{y \to x^-} (f(y)) = f(x)$ for all $x \in (0,1]$) is called piecewise $C^1$ on $[0,1]$ and we write $f \in PC^1([0,1])$, if the following properties hold: (i) for every $x \in [0,1)$ the limits $\lim_{y \to x^+} (f(y))$, $\lim_{h \to 0^+, y \to x^+} \left( h^{-1}(f(y+h) - f(y)) \right)$ exist and are finite, (ii) for every $x \in (0,1]$ the limit $\lim_{h \to 0^-} \left( h^{-1}(f(x+h) - f(x)) \right)$ exists and is finite, (iii) there exists a set $J \subset (0,1)$ of finite cardinality, where $\frac{df}{dx}(x) = \lim_{h \to 0^-} \left( h^{-1}(f(x+h) - f(x)) \right) = \lim_{h \to 0^+} \left( h^{-1}(f(x+h) - f(x)) \right)$ holds for $x \in (0,1) \setminus J$, and (iv) the mapping $(0,1) \setminus J \ni x \to \frac{df}{dx}(x) \in \Re$ is continuous. Notice that we require a piecewise $C^1$ function to be left-continuous but not continuous.

## 2. Stability Estimates for the Continuity Equation

We consider the continuity equation on a bounded domain, i.e., we consider the equation

$$\frac{\partial \rho}{\partial t}(t,x) + v(t,x)\frac{\partial \rho}{\partial x}(t,x) + \rho(t,x)\frac{\partial v}{\partial x}(t,x) = 0, \text{ for } (t,x) \in \Re_+ \times [0,1] \quad (2.1)$$

where

- $\rho$ is the state, required to be positive (i.e., $\rho(t,x) > 0$ for $(t,x) \in \Re_+ \times [0,1]$) and having a spatially uniform nominal equilibrium profile $\rho(x) \equiv \rho_s$, where $\rho_s > 0$ is a constant,
- $v$ is an input which has a spatially uniform nominal profile $v(x) \equiv v_s$, where $v_s > 0$ is a constant.

Equation (2.1) is accompanied by the boundary condition

$$\rho(t,0) = \rho_s \exp(b(t)), \text{ for } t \geq 0 \quad (2.2)$$

where $b$ is an additional input (the boundary disturbance).

For system (2.1), (2.2), we obtain the following result.

**Theorem 2.1:** *Consider the initial-boundary value problem (2.1), (2.2) with*

$$\rho(0,x) = \rho_0(x), \text{ for } x \in (0,1] \quad (2.3)$$

*where $\rho_s > 0$ is a constant, $\rho_0 \in C^1([0,1];(0,+\infty))$, $b \in C^1(\Re_+)$, and $v \in C^1(\Re_+ \times [0,1];(0,+\infty))$ is a positive function with $\rho_s \exp(b(0)) = \rho_0(0)$, $\frac{\partial v}{\partial x}(0,0) = -\dot{b}(0) - v(0,0)\frac{\rho_0'(0)}{\rho_0(0)}$. Then there exists a unique function $\rho \in C^1(\Re_+ \times [0,1];(0,+\infty))$ such that equations (2.1), (2.2), (2.3) hold. Furthermore, the function $\rho \in C^1(\Re_+ \times [0,1];(0,+\infty))$ satisfies the following estimates for all $t \geq 0$, $p \in (1,+\infty)$, $\mu > 0$ with $\mu > -p^{-1}v_{\max}(t) - v_{\min}(t)$:*

$$\left(\int_0^1 \left|\ln\left(\frac{\rho(t,x)}{\rho_s}\right)\right|^p ds\right)^{1/p} \leq \exp\left(p^{-1}v_{\max}(t)t\right) h\left(t - \frac{1}{v_{\min}(t)}\right)\left(\int_0^1 \left|\ln\left(\frac{\rho_0(x)}{\rho_s}\right)\right|^p ds\right)^{1/p}$$

$$+ \frac{1}{v_{\min}(t)}\exp\left(1 + \frac{\mu + p^{-1}v_{\max}(t)}{v_{\min}(t)}\right) \max_{\max\left(0, t - \frac{1}{v_{\min}(t)}\right) \leq s \leq t}\left(\left\|\frac{\partial v}{\partial x}[s]\right\|_p \exp(-\mu(t-s))\right) \quad (2.4)$$

$$+ \left(v_{\min}(t)\frac{\exp\left(\frac{p\mu}{v_{\min}(t)}\right) - 1}{p\mu}\right)^{1/p} \max_{\max\left(0, t - \frac{1}{v_{\min}(t)}\right) \leq s \leq t}\left(|b(s)|\exp(-\mu(t-s))\right)$$

$$\max_{0 \leq x \leq 1}\left(\left|\ln\left(\frac{\rho(t,x)}{\rho_s}\right)\right|\right) \leq h\left(t - \frac{1}{v_{\min}(t)}\right)\max_{0 \leq x \leq 1}\left(\left|\ln\left(\frac{\rho_0(x)}{\rho_s}\right)\right|\right)$$

$$+ \frac{1}{v_{\min}(t)}\exp\left(1 + \frac{\mu}{v_{\min}(t)}\right)\max_{\max\left(0, t - \frac{1}{v_{\min}(t)}\right) \leq s \leq t}\left(\left\|\frac{\partial v}{\partial x}[s]\right\|_\infty \exp(-\mu(t-s))\right) \quad (2.5)$$

$$+ \exp\left(\frac{\mu}{v_{\min}(t)}\right)\max_{\max\left(0, t - \frac{1}{v_{\min}(t)}\right) \leq s \leq t}\left(|b(s)|\exp(-\mu(t-s))\right)$$

*where* $h(s) := \begin{cases} 1 \text{ for } s < 0 \\ 0 \text{ for } s \geq 0 \end{cases}$ *and*

$$v_{\min}(t) = \min\{v(s,x) : s \in [0,t], x \in [0,1]\}$$
$$v_{\max}(t) = \max\left\{\frac{\partial v}{\partial x}(s,x) : s \in [0,t], x \in [0,1]\right\}$$
(2.6)

*for* $t \geq 0$.

**Remark 2.2:** **(a)** Estimates (2.4), (2.5) are stability estimates in special state norms. Due to the positivity of the state, the logarithmic norm of the state $\rho$ appears, i.e., we have $\left|\ln\left(\frac{\rho(t,x)}{\rho_s}\right)\right|$ instead of the usual $|\rho(t,x) - \rho_s|$ that appears in many stability estimates for linear PDEs. The logarithmic norm is a manifestation of the nonlinearity of system (2.1), (2.2) and the fact that the state space is not a linear space but rather a positive cone: the state space for system (2.1), (2.2) is the set $X = C^1([0,1];(0,+\infty))$.
**(b)** The set of allowable inputs is not a linear space: it is the set of all $b \in C^1(\Re_+)$, and $v \in C^1(\Re_+ \times [0,1];(0,+\infty))$ with $\rho_s \exp(b(0)) = \rho_0(0)$, $\frac{\partial v}{\partial x}(0,0) = -\dot{b}(0) - v(0,0)\frac{\rho_0'(0)}{\rho_0(0)}$. Again this fact is a manifestation of the nonlinearity of system (2.1), (2.2).
**(c)** Estimate (2.4) is not an ISS estimate, due to the fact that the overshoot coefficient bounded by $\exp\left(\frac{v_{\max}(t)}{pv_{\min}(t)}\right)$ (notice that $\exp(p^{-1}v_{\max}(t)t)h\left(t - \frac{1}{v_{\min}(t)}\right) \leq \exp\left(\frac{v_{\max}(t)}{pv_{\min}(t)}\right)$ for all $t \geq 0$) depends heavily on the input $v$. Estimates (2.4), (2.5) indicate that the gain coefficient of the boundary disturbance $b$ also depends on the input $v$.
**(d)** In the absence of disturbances, i.e., when $b(t) \equiv 0$ and $v(t,x) \equiv v_s > 0$ both estimates (2.4), (2.5) indicate finite-time stability. This is a well-known phenomenon to linear transport PDEs (see [11,12]).
**(e)** There is another feature of the stability estimates (2.4), (2.5) that should be noticed: the feature of finite-time memory. Indeed, only the input values in the time interval $\left[\max\left(0, t - \frac{1}{v_{\min}(t)}\right), t\right]$ affect the state at time $t \geq 0$ and this is manifested in estimates (2.4), (2.5) by the use of the operators $\max_{\max\left(0, t - \frac{1}{v_{\min}(t)}\right) \leq s \leq t}\left(\left\|\frac{\partial v}{\partial x}[s]\right\|_\infty \exp(-\mu(t-s))\right)$ and $\max_{\max\left(0, t - \frac{1}{v_{\min}(t)}\right) \leq s \leq t}\left(|b(s)|\exp(-\mu(t-s))\right)$ instead of the usual operators $\max_{0 \leq s \leq t}\left(\left\|\frac{\partial v}{\partial x}[s]\right\|_\infty \exp(-\mu(t-s))\right)$ and $\max_{0 \leq s \leq t}\left(|b(s)|\exp(-\mu(t-s))\right)$ that would appear in a standard ISS estimate.
**(f)** Applying the constant velocity $v(t,x) \equiv v_s$, the constant boundary disturbance $b(t) \equiv b$ and the constant initial condition $\rho_0(x) \equiv \rho_s \exp(b)$, one can show that the gain of the boundary disturbance has to be greater or equal to 1 for every $L^p(0,1)$ norm with $p \in (1, +\infty]$. On the other hand, using (2.4), (2.5) we obtain that the gain of this specific boundary disturbance has to be less than or equal to $\left(v_s \frac{\exp\left(\frac{p\mu}{v_s}\right) - 1}{p\mu}\right)^{1/p}$ for every $L^p(0,1)$ norm with $p \in (1, +\infty)$ and less than or equal to $\exp\left(\frac{\mu}{v_s}\right)$ for $p = +\infty$, where $\mu > 0$ is arbitrary. Since $\lim_{\mu \to 0^+} \left( v_s \dfrac{\exp\left(\dfrac{p\mu}{v_s}\right) - 1}{p\mu} \right)^{1/p} = 1$, it follows that the estimation of the gain of the boundary disturbance is optimal for this case.

**(g)** The time-invariant velocities $v_1(x) = 1 + (\theta - 1)x$ and $v_2(x) = \theta + (1 - \theta)x$ with $\theta \in (0,1)$, have equal minimum and maximum values. Moreover, $\left\| \dfrac{\partial v_i}{\partial x}[t] \right\|_p = 1 - \theta$ for all $p \in (1, +\infty]$ for $i = 1, 2$. Applying these velocities, the constant boundary disturbance $b(t) \equiv 0$ and the initial conditions $\rho_{0,i}(x) \equiv \rho_s \dfrac{v_i(0)}{v_i(x)}$ for $i = 1, 2$, one can show that the corresponding gains of $\dfrac{\partial v}{\partial x}$ have to be greater than or equal to

$$\gamma_1 = \dfrac{1}{1-\theta} \int_0^1 \left(-\ln(1 + (\theta - 1)x)\right)^p dx \text{ and } \gamma_2 = \dfrac{1}{1-\theta} \int_0^1 \left(\ln(1 + (\theta^{-1} - 1)x)\right)^p dx,$$

respectively, for every $L^p(0,1)$ norm with $p \in (1, +\infty)$. Since $\dfrac{1}{1 + (\theta - 1)x} < 1 + (\theta^{-1} - 1)x$ for all $x \in (0,1)$, it follows that the gain of $v_2(x) = \theta + (1 - \theta)x$ is strictly greater than the gain of $v_1(x) = 1 + (\theta - 1)x$. Therefore, velocities that are increasing with respect to $x$ (i.e., convection that speeds up downstream) add a greater bias to the solution profile than velocities that are decreasing with respect to $x$ (i.e., convection that slows down downstream). This is also apparent from the estimation of the gain of $\dfrac{\partial v}{\partial x}$ from (2.4) for every $L^p(0,1)$ norm with $p \in (1, +\infty)$: the gain estimate $\dfrac{1}{v_{\min}(t)} \exp\left(1 + \dfrac{\mu + p^{-1} v_{\max}(t)}{v_{\min}(t)}\right)$ depends on $v_{\max}(t) = \max\left\{ \dfrac{\partial v}{\partial x}(s, x) : s \in [0,t], x \in [0,1] \right\}$ and indicates that velocities that are increasing with respect to $x$ (i.e., convection that speeds up downstream with $v_{\max}(t) \geq 0$) add a greater bias to the solution profile than velocities that are decreasing with respect to $x$ (i.e., convection that slows down downstream for which $v_{\max}(t) \leq 0$).

**(h)** When the inputs are constant in time, i.e., $b(t) \equiv b$ and $v(t,x) \equiv v(x)$, then the equilibrium profiles of system (2.1), (2.2) are given by the equation $\rho(x) = \rho_s \exp(b) \dfrac{v(0)}{v(x)}$ for $x \in [0,1]$. Consequently, it becomes clear that the velocity $v$ acts as an "equilibrium-shaping functional parameter". For instance, if $v(x)$ is monotonically increasing, namely, if the convection speeds up downstream, the density equilibrium profile $\rho(x)$ decreases — and vice versa — and such non constant equilibrium profiles are finite-time stable (in logarithmic norm).

The proof of Theorem 2.1 relies on the following result which has its own interest.

**Theorem 2.3:** *Consider the initial-boundary value problem*

$$\dfrac{\partial w}{\partial t}(t,x) + v(t,x) \dfrac{\partial w}{\partial x}(t,x) = a(t,x) w(t,x) + f(t,x) \tag{2.7}$$

$$w(0,x) = \varphi(x), \text{ for } x \in (0,1] \tag{2.8}$$

$$w(t,0) = b(t), \text{ for } t \geq 0 \tag{2.9}$$

where $\varphi \in PC^1([0,1])$, $b \in C^1(\mathfrak{R}_+)$, $f \in C^1(\mathfrak{R}_+ \times [0,1])$, $a \in C^1(\mathfrak{R}_+ \times [0,1])$ and $v \in C^1(\mathfrak{R}_+ \times [0,1];(0,+\infty))$ is a positive function. Let $\xi_i \in [0,1)$ ($i = 0,...,N$) with $\xi_0 = 0$ be the points for which $\varphi \in C^1([0,1] \setminus \{\xi_0,...,\xi_N\})$. Let $r_i$ ($i = 0,...,N$) be the solutions of the initial value problems $\dot{r}_i(t) = v(t, r_i(t))$ with $r_i(0) = \xi_i$ and if there exists $T_i > 0$ with $r_i(T_i) = 1$ then define $r_i(t) = 1$ for all $t > T_i$. Then there exists a unique function $w: \mathfrak{R}_+ \times [0,1] \to \mathfrak{R}$ of class $C^1(\mathfrak{R}_+ \times [0,1] \setminus \Omega)$, where $\Omega = \bigcup_{i=0,...,N} \{(t, r_i(t)) : t \geq 0, r_i(t) < 1\}$ with $w[t] \in PC^1([0,1])$ for all $t \geq 0$, such that (2.7) holds for all $(t,x) \in \mathfrak{R}_+ \times [0,1] \setminus \Omega$ and equations (2.8), (2.9) hold. Furthermore, the function $w: \mathfrak{R}_+ \times [0,1] \to \mathfrak{R}$ satisfies the following estimates for all $t \geq 0$, $p \in (1,+\infty)$, $\mu \geq 0$ with $\mu > -A(t)$, $\mu > -p^{-1} v_{\max}(t) - A(t) - v_{\min}(t)$:

$$\|w[t]\|_p \leq \exp\left((A(t) + p^{-1} v_{\max}(t))t\right) h\left(t - \frac{1}{v_{\min}(t)}\right) \|\varphi\|_p$$
$$+ \frac{1}{v_{\min}(t)} \exp\left(1 + \frac{\mu + p^{-1} v_{\max}(t) + A(t)}{v_{\min}(t)}\right) \max_{\max\left(0, t - \frac{1}{v_{\min}(t)}\right) \leq s \leq t} \left(\|f[s]\|_p \exp(-\mu(t-s))\right) \quad (2.10)$$
$$+ \left(v_{\min}(t) \frac{\exp\left(\frac{p(\mu + A(t))}{v_{\min}(t)}\right) - 1}{p(\mu + A(t))}\right)^{1/p} \max_{\max\left(0, t - \frac{1}{v_{\min}(t)}\right) \leq s \leq t} \left(|b(s)| \exp(-\mu(t-s))\right)$$

$$\|w[t]\|_\infty \leq \exp(t A(t)) h\left(t - \frac{1}{v_{\min}(t)}\right) \|\varphi\|_\infty$$
$$+ \frac{1}{v_{\min}(t)} \exp\left(1 + \frac{\mu + A(t)}{v_{\min}(t)}\right) \max_{\max\left(0, t - \frac{1}{v_{\min}(t)}\right) \leq s \leq t} \left(\|f[s]\|_\infty \exp(-\mu(t-s))\right) \quad (2.11)$$
$$+ \exp\left(\frac{\mu + A(t)}{v_{\min}(t)}\right) \max_{\max\left(0, t - \frac{1}{v_{\min}(t)}\right) \leq s \leq t} \left(|b(s)| \exp(-\mu(t-s))\right)$$

where $h(s) := \begin{cases} 1 \text{ for } s < 0 \\ 0 \text{ for } s \geq 0 \end{cases}$ and

$$v_{\min}(t) = \min\{v(s,x) : s \in [0,t], x \in [0,1]\}$$
$$v_{\max}(t) = \max\left\{\frac{\partial v}{\partial x}(s,x) : s \in [0,t], x \in [0,1]\right\} \quad (2.12)$$
$$A(t) = \max\{a(s,x) : s \in [0,t], x \in [0,1]\}$$

for all $t \geq 0$. Moreover, if $b(0) = \varphi(0)$ and $\varphi \in C^0([0,1])$ then $w \in C^0(\mathfrak{R}_+ \times [0,1])$. Finally, if $\varphi \in C^1([0,1])$, $b(0) = \varphi(0)$, $\dot{b}(0) + v(0,0) \varphi'(0) = a(0,0) b(0) + f(0,0)$ then $w \in C^1(\mathfrak{R}_+ \times [0,1])$.

The proof of Theorem 2.3 is provided in Section 4 and is based on a combination of different methodologies:

- The exploitation of the superposition principle for the initial-value problem (2.7), (2.8), (2.9): the solution of (2.7), (2.8), (2.9) can be written as the sum of three functions: the solution of (2.7), (2.8), (2.9) with zero inputs $f, b$ and initial condition the given function $\varphi \in PC^1([0,1])$, the solution of (2.7), (2.8), (2.9) with zero initial condition, zero distributed input $f$ and boundary input the given function $b \in C^1(\Re_+)$, and the solution of (2.7), (2.8), (2.9) with zero initial condition, zero boundary input $f$ and distributed input the given function $f \in C^1(\Re_+ \times [0,1])$.

- The norms of the first two components of the solution are estimated by using the exact formulas of the solution on the characteristic curves of the PDE (2.7).

- The norm of the third component of the solution is estimated by using a Lyapunov analysis.

## 3. Feedback Control of Manufacturing Systems

Manufacturing systems with a high volume and a large number of consecutive production steps (which typically number in the many hundreds) are often modelled by non-local PDEs of the form (see [7,8,9,23]):

$$\frac{\partial \rho}{\partial t}(t,x) + \lambda(W(t))\frac{\partial \rho}{\partial x}(t,x) = 0 \text{, for } (t,x) \in \Re_+ \times [0,1] \tag{3.1}$$

$$W(t) = \int_0^1 \rho(t,x)dx \text{, for } t \geq 0 \tag{3.2}$$

where

- $\rho(t,x)$ is the density of the processed material at time $t \geq 0$ and stage $x \in [0,1]$, required to be positive (i.e., $\rho(t,x) > 0$ for $(t,x) \in \Re_+ \times [0,1]$) and having a spatially uniform equilibrium profile $\rho(x) \equiv \rho_s$, where $\rho_s > 0$ is a constant (set point), and
- $\lambda \in C^1((0,+\infty);(0,+\infty))$ is a nonlinear function that determines the production speed.

The model is accompanied by the influx boundary condition

$$\rho(t,0)\lambda(W(t)) = u(t) \text{, for } t \geq 0 \tag{3.3}$$

where $u(t) \in (0,+\infty)$ is the control input (the process influx rate). Existence, uniqueness and related control problems for systems of the form (3.1), (3.2), (3.3) were studied in [7,8,9]. Here we want to address the feedback stabilization problem of the spatially uniform equilibrium profile $\rho(x) \equiv \rho_s$ under the feedback control law

$$u(t) = \rho_s \lambda(W(t))\exp(b(t)) \text{, for } t \geq 0 \tag{3.4}$$

where $b(t) \in \Re$ represents an uncertainty. The motivation for the feedback law (3.4) comes from the fact that in the absence of uncertainties the feedback law (3.4) combined with the boundary condition (3.3) gives the boundary condition $\rho(t,0) = \rho_s$ for $t \geq 0$, which guarantees finite-time stability for the linearization of (3.1), (3.2). Moreover, the implementation of the feedback law (3.4) relies on the measurement of the total load in the production line $W(t)$ which is a quantity that can be measured relatively easily. Finally, notice that the feedback law (3.4) guarantees that $u(t)$ is positive and for bounded disturbances as well as bounded functions $\lambda \in C^1((0,+\infty);(0,+\infty))$ (which is usually the case for manufacturing systems) the control input $u(t)$ is bounded from above by a constant independent of the initial condition (bounded feedback).

Using Theorem 2.1, we are in a position to obtain the following result.

**Theorem 3.1:** *Consider the initial-boundary value problem (3.1), (3.2), (3.3), (3.4) with*

$$\rho(0,x) = \rho_0(x), \text{ for } x \in (0,1] \tag{3.5}$$

*where $\rho_s > 0$ is a constant, $\rho_0 \in C^1([0,1];(0,+\infty))$ and $b \in C^1(\Re_+)$ is a bounded function with $\rho_s \exp(b(0)) = \rho_0(0)$, $\dot{b}(0) = -\lambda \left( \int_0^1 \rho_0(x) dx \right) \frac{\rho_0'(0)}{\rho_0(0)}$. Then there exists a unique function $\rho \in C^1(\Re_+ \times [0,1];(0,+\infty))$ such that equations (3.1), (3.2), (3.3), (3.4), (3.5) hold. Furthermore, the function $\rho \in C^1(\Re_+ \times [0,1];(0,+\infty))$ satisfies the following estimates for all $t \geq 0$, $p \in (1,+\infty)$, $\mu > 0$:*

$$\left( \int_0^1 \left| \ln\left( \frac{\rho(t,x)}{\rho_s} \right) \right|^p ds \right)^{1/p} \leq h(t-r) \left( \int_0^1 \left| \ln\left( \frac{\rho_0(x)}{\rho_s} \right) \right|^p ds \right)^{1/p}$$
$$+ \left( \frac{\exp(p\mu r) - 1}{p\mu r} \right)^{1/p} \max_{\max(0,t-r) \leq s \leq t} \left( |b(s)| \exp(-\mu(t-s)) \right) \tag{3.6}$$

$$\max_{0 \leq x \leq 1} \left( \left| \ln\left( \frac{\rho(t,x)}{\rho_s} \right) \right| \right) \leq h(t-r) \max_{0 \leq x \leq 1} \left( \left| \ln\left( \frac{\rho_0(x)}{\rho_s} \right) \right| \right)$$
$$+ \exp(\mu r) \max_{\max(0,t-r) \leq s \leq t} \left( |b(s)| \exp(-\mu(t-s)) \right) \tag{3.7}$$

*where $h(s) := \begin{cases} 1 \text{ for } s < 0 \\ 0 \text{ for } s \geq 0 \end{cases}$ and*

$$r := \frac{1}{\min \left\{ \lambda(s) : \min\left( \min_{0 \leq x \leq 1} (\rho_0(x)), \rho_s \exp\left( \inf_{t \geq 0} (b(t)) \right) \right) \leq s \leq \max\left( \max_{0 \leq x \leq 1} (\rho_0(x)), \rho_s \exp\left( \sup_{t \geq 0} (b(t)) \right) \right) \right\}} \tag{3.8}$$

Estimates (3.6), (3.7) guarantee finite-time stability in the absence of uncertainties. However, notice that the terminal time $r$ given by (3.8) depends on the initial condition and the boundary disturbance. Therefore, for certain initial conditions or for large boundary disturbances it may happen that the terminal time $r$ is unacceptably large. For example, when $\lambda(W) = (1+W)^{-1}$ and $\rho_0 \in C^1([0,1];(0,+\infty))$ is any function (its monotonicity does not play any role) then (3.8) gives $r := 1 + \max\left( \|\rho_0\|_\infty, \rho_s \exp\left( \sup_{t \geq 0} (b(t)) \right) \right)$. Consequently, if the initial density is large for some $x \in [0,1]$ or if the boundary disturbance is large, then the terminal time $r$ will be large. This is a possible disadvantage of the feedback law (3.4) and we do not know whether it is possible to achieve smaller terminal times (see also the discussion in Section 4 of [7] for the case $\lambda(W) = (1+W)^{-1}$). Estimates (3.6), (3.7) guarantee robustness with respect to the boundary uncertainty $b$. However, again when the terminal time $r$ is very large then the gain of the uncertainty becomes very large.

# 4. Proofs of Main Results

The proof of Theorem 2.1 is simply an application of Theorem 2.3 and the use of the transformation $w(t,x) = \ln\left(\frac{\rho(t,x)}{\rho_s}\right)$ (or its inverse $\rho(t,x) = \rho_s \exp(w(t,x))$). Therefore, we next focus on the proof of Theorem 2.3.

**Proof of Theorem 2.3:** Extending $v \in C^1(\Re_+ \times [0,1];(0,+\infty))$ so that $v \in C^1(\Re^2)$ we can define for all $t_0 \geq 0$, $x_0 \in [0,1)$ the mapping $X(s;t_0,x_0) \in [0,1]$ as the unique solution of the initial-value problem

$$\frac{dX}{ds}(s) = v(t_0 + s, X(s)), \quad X(t_0) = x_0 \tag{4.1}$$

Clearly, $X(s;t_0,x_0) \in [0,1]$ is defined for $s \in [0, s_{\max})$, where $s_{\max} \in (0,+\infty]$ is the maximal existence time of the solution. Since $v(t,x) > 0$ for all $t \geq 0$, $x \in [0,1]$, it follows that a finite $s_{\max}$ implies that $\lim_{s \to s_{\max}^-}(X(s;t_0,x_0)) = X(s_{\max};t_0,x_0) = 1$. Notice that the mapping $X(s;t_0,x_0)$ is increasing with respect to $s \in [0, s_{\max})$ and $x_0 \in [0,1)$ and satisfies the following equations for all $t_0 \geq 0$, $x_0 \in [0,1)$ and $s \in [0, s_{\max})$:

$$\frac{\partial X}{\partial x_0}(s;t_0,x_0) = \exp\left(\int_0^s \frac{\partial v}{\partial x}(t_0 + l, X(l;t_0,x_0))dl\right) > 0$$

$$\frac{\partial X}{\partial t_0}(s;t_0,x_0) = v(t_0 + s, X(s;t_0,x_0)) - v(t_0, x_0)\exp\left(\int_0^s \frac{\partial v}{\partial x}(t_0 + l, X(l;t_0,x_0))dl\right) \tag{4.2}$$

It follows from (4.1) and (4.2) (which imply that $\frac{\partial}{\partial t_0}X(t - t_0;t_0,0) = -v(t_0,0)\exp\left(\int_0^{t-t_0} \frac{\partial v}{\partial x}(t_0 + l, X(l;t_0,0))dl\right) < 0$) that

- for every $t \geq 0$, $x \in [0,1)$ with $x > r_0(t) = X(t;0,0)$ the equation $X(t;0,x_0) = x$ will be uniquely solvable with respect to $x_0 \in [0,1)$, and
- for every $t \geq 0$, $x \in [0,1)$ with $x \leq r_0(t) = X(t;0,0)$ the equation $X(t - t_0;t_0,0) = x$ will be uniquely solvable with respect to $t_0 \geq 0$.

Let $x_0(t,x) \in [0,1]$ and $t_0(t,x) \geq 0$ be the solutions of the above equations. By virtue of the implicit function theorem they are both $C^1$ on their domains. We define:

$$w(t,x) := \exp\left(\int_0^t a(s, X(s;0,x_0(t,x)))ds\right)\varphi(x_0(t,x))$$

$$+ \int_0^t \exp\left(\int_\tau^t a(s, X(s;0,x_0(t,x)))ds\right) f(\tau, X(\tau;0,x_0(t,x)))d\tau$$

for $t \geq 0$, $x \in [0,1)$ with $x > r_0(t) = X(t;0,0)$ \hfill (4.3)

and

$$w(t,x) := \exp\left(\int_{t_0(t,x)}^t a(s, X(s-t_0(t,x);t_0(t,x),0))ds\right)b(t_0(t,x))$$

$$+ \int_{t_0(t,x)}^t \exp\left(\int_\tau^t a(s, X(s-t_0(t,x);t_0(t,x),0))ds\right) f(\tau, X(\tau-t_0(t,x);t_0(t,x),0))d\tau$$

for $t \geq 0$, $x \in [0,1)$ with $x \leq r_0(t) = X(t;0,0)$ \hfill (4.4)

Definitions (4.3), (4.4) guarantee that $w: \Re_+ \times [0,1] \to \Re$ is a function of class $C^1(\Re_+ \times [0,1] \setminus \Omega)$, where $\Omega = \cup_{i=0,...,N} \{(t, r_i(t)): t \geq 0, r_i(t) < 1\}$ with $w[t] \in PC^1([0,1])$ for all $t \geq 0$, such that (2.7) holds for all $(t,x) \in \Re_+ \times [0,1] \setminus \Omega$ and equations (2.8), (2.9) hold. Moreover, definitions (4.3), (4.4) guarantee that if $b(0) = \varphi(0)$ and $\varphi \in C^0([0,1])$ then $w \in C^0(\Re_+ \times [0,1])$ and if $\varphi \in C^1([0,1])$, $b(0) = \varphi(0)$, $\dot{b}(0) + v(0,0)\varphi'(0) = a(0,0)b(0) + f(0,0)$ then $w \in C^1(\Re_+ \times [0,1])$.

Uniqueness follows from a contradiction argument. Suppose that there exist two functions $w, \tilde{w} \in C^1(\Re_+; L^2(0,1))$ with $w[t], \tilde{w}[t] \in PC^1([0,1])$ for $t \geq 0$, that satisfy (2.7) (in the $L^2(0,1)$ sense) and (2.8), (2.9). It then follows that the function $e = w - \tilde{w}$ satisfies the following equations:

$$\frac{\partial e}{\partial t}(t,x) + v(t,x)\frac{\partial e}{\partial x}(t,x) = a(t,x)e(t,x), \text{ for } t \geq 0 \tag{4.5}$$

$$e(0,x) = 0, \text{ for } x \in (0,1] \tag{4.6}$$

$$e(t,0) = 0, \text{ for } t \geq 0 \tag{4.7}$$

Using the functional $V(t) = \int_0^1 e^2(t,x)dx$ on $[0,T]$ for arbitrary $T > 0$, we have by virtue of (2.12), (4.5) and (4.7) for every $t \in [0,T]$:

$$\dot{V}(t) = 2\int_0^1 e(t,x)\frac{\partial e}{\partial t}(t,x)dx = -2\int_0^1 v(t,x)e(t,x)\frac{\partial e}{\partial x}(t,x)dx + 2\int_0^1 a(t,x)e^2(t,x)dx$$

$$= -\int_0^1 v(t,x)\frac{\partial}{\partial x}\left(e^2(t,x)\right)dx + 2A(T)V(t)$$

$$= -v(t,1)e^2(t,1) + \int_0^1 e^2(t,x)\frac{\partial v}{\partial x}(t,x)dx + 2A(T)V(t)$$

$$\leq (v_{\max}(T) + 2A(T))V(t)$$

Gronwall's lemma implies that $V(t) \leq \exp\left((v_{\max}(T) + 2A(T))t\right)V(0)$ for all $t \in [0,T]$ and consequently (using (4.6)), we get $V(t) = 0$ for all $t \in [0,T]$. This equality in conjunction with the fact that $w[t], \tilde{w}[t] \in PC^1([0,1])$ for $t \geq 0$, implies $w \equiv \tilde{w}$.

Using (4.3), (4.4), we next notice that

$$w[t] = w_1[t] + w_2[t] + w_3[t], \text{ for } t \geq 0 \tag{4.8}$$

where

$$w_1(t,x) := 0, \quad w_2(t,x) := \exp\left(\int_0^t a(s, X(s;0,x_0(t,x)))ds\right)\varphi(x_0(t,x)),$$

$$w_3(t,x) := \int_0^t \exp\left(\int_\tau^t a(s, X(s;0,x_0(t,x)))ds\right) f(\tau, X(\tau;0,x_0(t,x)))d\tau$$

$$\text{for } t \geq 0, \ x \in [0,1) \text{ with } x > r_0(t) = X(t;0,0) \tag{4.9}$$

and

$$w_1(t,x) := \exp\left(\int_{t_0(t,x)}^t a(s, X(s-t_0(t,x); t_0(t,x), 0))ds\right) b(t_0(t,x)), \quad w_2(t,x) := 0,$$

$$w_3(t,x) := \int_{t_0(t,x)}^t \exp\left(\int_\tau^t a(s, X(s-t_0(t,x); t_0(t,x), 0))ds\right) f(\tau, X(\tau-t_0(t,x); t_0(t,x), 0))d\tau$$

$$\text{for } t \geq 0, \ x \in [0,1) \text{ with } x \leq r_0(t) = X(t;0,0) \tag{4.10}$$

First we estimate the $L^p(0,1)$ norm of $w_1$ with $p \in [1,+\infty)$. Notice that by virtue of (2.12), (4.1) and since $t_0(t,x) \geq 0$ solves the equation $X(t-t_0;t_0,0) = x$, we get for all $t \geq 0$, $x \in [0,1)$ with $x \leq r_0(t) = X(t;0,0)$:

$$\max\left(0, t - \frac{x}{v_{\min}(t)}\right) \leq t_0(t,x) \leq t \tag{4.11}$$

Therefore, we get from (4.10), (2.12), (4.11) for every $\mu \geq \max(0, -A(t))$, $p \in [1,+\infty)$ and $t \geq 0$, $x \in [0,1)$ with $x \leq r_0(t) = X(t;0,0)$:

$$\begin{aligned}
|w_1(t,x)|^p &= \exp\left(p \int_{t_0(t,x)}^{t} a(s, X(s-t_0(t,x); t_0(t,x), 0))ds\right)|b(t_0)|^p \\
&\leq \exp\left(pA(t)(t-t_0(t,x))\right)|b(t_0(t,x))|^p \\
&\leq \exp\left(p(\mu+A(t))(t-t_0(t,x))\right) \max_{t_0(t,x) \leq s \leq t}\left(|b(s)|^p \exp(-p\mu(t-s))\right) \\
&\leq \exp\left(\frac{p(\mu+A(t))x}{v_{\min}(t)}\right) \max_{\max\left(0, t - \frac{1}{v_{\min}(t)}\right) \leq s \leq t}\left(|b(s)|^p \exp(-p\mu(t-s))\right)
\end{aligned} \tag{4.12}$$

Using (4.9) and (4.12), we obtain for every $\mu \geq 0$ with $\mu > -A(t)$, $p \in [1,+\infty)$ and $t \geq 0$:

$$\|w_1[t]\|_p \leq \left(v_{\min}(t) \frac{\exp\left(\frac{p(\mu+A(t))}{v_{\min}(t)}\right) - 1}{p(\mu+A(t))}\right)^{1/p} \max_{\max\left(0, t - \frac{1}{v_{\min}(t)}\right) \leq s \leq t}\left(|b(s)| \exp(-\mu(t-s))\right) \tag{4.13}$$

Next we estimate the $L^p(0,1)$ norm of $w_2$ with $p \in [1,+\infty)$. For all $t \geq 0$ with $r_0(t) < 1$ we get from (4.9) and (2.12):

$$|w_2(t,x)| \leq \exp(A(t)t)|\varphi(x_0(t,x))|, \text{ for } x \in [0,1) \text{ with } x > r_0(t) \tag{4.14}$$

Using the fact that $\frac{\partial x_0}{\partial x}(t,x) = \exp\left(-\int_0^t \frac{\partial v}{\partial x}(s, X(s;0,x_0(t,x)))ds\right)$ (a consequence of (4.1), (4.2) and the fact that $X(t;0,x_0(t,x)) = x$), we get from (4.9), (4.10), (2.12), (4.14) and the substitution $\xi = x_0(t,x)$ (allowable since $\varphi \in PC^1([0,1])$) for all $p \in [1,+\infty)$ and $t \geq 0$ with $r_0(t) < 1$:

$$\begin{aligned}
\|w_2[t]\|_p &\leq \exp(A(t)t) \left(\int_{r_0(t)}^{1} |\varphi(x_0(t,x))|^p dx\right)^{1/p} \\
&\leq \exp(A(t)t) \left(\int_0^1 |\varphi(\xi)|^p \exp\left(\int_0^t \frac{\partial v}{\partial x}(s, X(s;0,\xi))ds\right)d\xi\right)^{1/p} \\
&\leq \exp\left((A(t) + p^{-1}v_{\max}(t))t\right)\|\varphi\|_p
\end{aligned} \tag{4.15}$$

Notice that the existence of $t \geq 0$, $x \in [0,1)$ with $x > r_0(t) = X(t;0,0)$ in conjunction with (2.12) (which gives $r_0(t) \geq tv_{\min}(t)$) implies that $t < \frac{1}{v_{\min}(t)}$. Consequently, when $t \geq \frac{1}{v_{\min}(t)}$ then it follows from (4.10) that $\|w_2[t]\|_p = 0$ for all $p \in [1,+\infty)$. Thus we obtain from (4.15) for all $t \geq 0$ and $p \in [1,+\infty)$:

$$\|w_2[t]\|_p \leq \exp\left((A(t) + p^{-1}v_{\max}(t))t\right) h\left(t - \frac{1}{v_{\min}(t)}\right) \|\varphi\|_p \qquad (4.16)$$

where $h(s) := \begin{cases} 1 & \text{for } s < 0 \\ 0 & \text{for } s \geq 0 \end{cases}$.

Next we estimate the $L^p(0,1)$ norm of $w_3$ with $p \in [1,+\infty)$. Notice that $w_3$ is the solution of (2.7), (2.8), (2.9) with $\varphi \equiv 0$ and $b \equiv 0$. Let $\sigma > 0$ be a constant (to be selected) and define for $p \in (1,+\infty)$, $t \geq 0$:

$$V(t) = \int_0^1 \exp(-\sigma x) |w_3(t,x)|^p \, dx \qquad (4.17)$$

Using the fact that $w_3$ is the solution of (2.7), (2.8), (2.9) with $\varphi \equiv 0$ and $b \equiv 0$, we get from (4.17) for all $t \geq 0$ and $s \in [0,t]$:

$$\dot{V}(s) = -p\int_0^1 \exp(-\sigma x) \operatorname{sgn}(w_3(s,x)) |w_3(s,x)|^{p-1} v(s,x) \frac{\partial w_3}{\partial x}(s,x) dx$$
$$+ p\int_0^1 \exp(-\sigma x) a(s,x) |w_3(s,x)|^p \, dx + p\int_0^1 \exp(-\sigma x) \operatorname{sgn}(w_3(s,x)) |w_3(s,x)|^{p-1} f(s,x) dx \qquad (4.18)$$

Integrating by parts, we obtain from (2.12), (4.17), (4.18) for all $t \geq 0$ and $s \in [0,t]$:

$$\dot{V}(s) \leq -\exp(-\sigma) v(s,1) |w_3(s,1)|^p + \int_0^1 \exp(-\sigma x) |w_3(s,x)|^p \frac{\partial v}{\partial x}(s,x) dx$$
$$- \sigma \int_0^1 \exp(-\sigma x) |w_3(s,x)|^p v(t,x) dx + pA(t)V(s) + p\int_0^1 \exp(-\sigma x) |w_3(s,x)|^{p-1} |f(s,x)| dx \qquad (4.19)$$

Using the inequality $|w(s,x)|^{p-1} |f(s,x)| \leq \frac{p-1}{p} \varepsilon^{\frac{p}{p-1}} |w(s,x)|^p + \frac{1}{p} \varepsilon^{-p} |f(s,x)|^p$ which holds for all $\varepsilon > 0$, we obtain from (2.12), (4.17), (4.19) for all $t \geq 0$, $\varepsilon > 0$ and $s \in [0,t]$:

$$\dot{V}(s) \leq -\left(\sigma v_{\min}(t) - v_{\max}(t) - pA(t) - (p-1)\varepsilon^{\frac{p}{p-1}}\right) V(s) + \varepsilon^{-p} \|f[s]\|_p^p \qquad (4.20)$$

Using Lemma 2.12 in [16] in conjunction with (4.20) and using the fact that $V(0) = 0$ (a consequence of (4.17) and the fact that $w_3$ is the solution of (2.7), (2.8), (2.9) with $\varphi \equiv 0$ and $b \equiv 0$), we obtain for all $t \geq 0$, $\varepsilon > 0$ and $\mu \geq 0$:

$$V(t) \leq \varepsilon^{-p} \int_0^t \exp\left(-(t-s)\left(\sigma v_{\min}(t) - v_{\max}(t) - pA(t) - (p-1)\varepsilon^{\frac{p}{p-1}}\right)\right) \|f[s]\|_p^p \, ds$$
$$\leq \varepsilon^{-p} \int_0^t \exp\left(-(t-s)\left(\sigma v_{\min}(t) - p\mu - v_{\max}(t) - pA(t) - (p-1)\varepsilon^{\frac{p}{p-1}}\right)\right) ds \max_{0 \leq s \leq t} \left(\|f[s]\|_p^p \exp(-p\mu(t-s))\right) \qquad (4.21)$$

Exploiting (4.21) in conjunction with (4.17), we get for all $t \geq 0$, $\varepsilon > 0$ and $\mu \geq 0$:

$$\|w_3[t]\|_p \leq K(t) \max_{0 \leq s \leq t} \left( \|f[s]\|_p \exp(-\mu(t-s)) \right) \tag{4.22}$$

where

$$K(t) := \exp\left(\frac{\sigma}{p}\right) \varepsilon^{-1} \left( \int_0^t \exp\left( -(t-s)\left( \sigma v_{\min}(t) - p\mu - v_{\max}(t) - pA(t) - (p-1)\varepsilon^{\frac{p}{p-1}} \right) \right) ds \right)^{1/p} \tag{4.23}$$

Selecting

$$\varepsilon = \left( \sigma p^{-1} v_{\min}(t) - \mu - p^{-1} v_{\max}(t) - A(t) \right)^{\frac{p-1}{p}}$$
$$\sigma = \frac{p\mu + v_{\max}(t) + pA(t) + pv_{\min}(t)}{v_{\min}(t)} \tag{4.24}$$

we obtain from (4.23) for all $t \geq 0$ and $\mu \geq 0$ with $\mu > -p^{-1}v_{\max}(t) - A(t) - v_{\min}(t)$:

$$\|w_3[t]\|_p \leq \frac{1}{v_{\min}(t)} \exp\left( 1 + \frac{\mu + p^{-1}v_{\max}(t) + A(t)}{v_{\min}(t)} \right) \max_{0 \leq s \leq t} \left( \|f[s]\|_p \exp(-\mu(t-s)) \right) \tag{4.25}$$

Noticing that $w_3[t]$ depends only on $f[s]$ with $s \in \left[ \max\left(0, t - \frac{1}{v_{\min}(t)}\right), t \right]$ and using a standard causality argument (when $t > \frac{1}{v_{\min}(t)}$) we obtain from (4.25) for all $t \geq 0$ and $\mu \geq 0$ with $\mu > -p^{-1}v_{\max}(t) - A(t) - v_{\min}(t)$:

$$\|w_3[t]\|_p \leq \frac{1}{v_{\min}(t)} \exp\left( 1 + \frac{\mu + p^{-1}v_{\max}(t) + A(t)}{v_{\min}(t)} \right) \max_{\max\left(0, t - \frac{1}{v_{\min}(t)}\right) \leq s \leq t} \left( \|f[s]\|_p \exp(-\mu(t-s)) \right) \tag{4.26}$$

Using (4.8), (4.13), (4.16), (4.26) and the triangle inequality, we obtain estimate (2.10). Estimate (2.11) is obtained by letting $p \to +\infty$ and by using the facts $\lim_{p \to +\infty}\left(\|w[t]\|_p\right) = \|w[t]\|_\infty$, $\|\varphi\|_p \leq \|\varphi\|_\infty$ for all $p \in [1, +\infty)$.

The proof is complete. ◁

Theorem 3.1 is a consequence of Theorem 2.1 and the following existence/uniqueness result which also provides a useful estimate.

**Proposition 4.1:** *Consider the initial-boundary value problem (3.1), (3.2), (3.5) with*

$$\rho(t, 0) = \rho_s \exp(b(t)), \text{ for } t \geq 0 \tag{4.27}$$

*where $\rho_s > 0$ is a constant, $\rho_0 \in C^1([0,1];(0,+\infty))$ and $b \in C^1(\Re_+)$ is a bounded function with $\rho_s \exp(b(0)) = \rho_0(0)$, $\dot{b}(0) = -\lambda\left(\int_0^1 \rho_0(x)dx\right)\frac{\rho_0'(0)}{\rho_0(0)}$. Then there exists a unique function $\rho \in C^1(\Re_+ \times [0,1];(0,+\infty))$ such that equations (3.1), (3.2), (3.5), (4.27) hold. Furthermore, the function $\rho \in C^1(\Re_+ \times [0,1];(0,+\infty))$ satisfies the following estimates for all $t \geq 0$, $x \in [0,1]$:*

$$\min\left( \min_{0 \leq y \leq 1}(\rho_0(y)), \rho_s \exp\left(\inf_{s \geq 0}(b(s))\right) \right) \leq \rho(t,x) \leq \max\left( \max_{0 \leq y \leq 1}(\rho_0(y)), \rho_s \exp\left(\sup_{s \geq 0}(b(s))\right) \right) \tag{4.28}$$

**Proof:** Define:

$$\tilde{b}(t) := \rho_s \exp(b(t)), \text{ for } t \geq 0 \tag{4.29}$$

$$\begin{aligned} v_{min} &= \min\{\lambda(s): \rho_{min} \leq s \leq \rho_{max}\} \\ v_{max} &= \max\{\lambda(s): \rho_{min} \leq s \leq \rho_{max}\} \end{aligned} \tag{4.30}$$

$$\begin{aligned} \rho_{min} &= \min\left(\min_{0 \leq y \leq 1}(\rho_0(y)), \inf_{s \geq 0}(\tilde{b}(s))\right) \\ \rho_{max} &= \max\left(\max_{0 \leq y \leq 1}(\rho_0(y)), \sup_{s \geq 0}(\tilde{b}(s))\right) \end{aligned} \tag{4.31}$$

$$L_\lambda = \max\{|\lambda'(s)|: \rho_{min} \leq s \leq \rho_{max}\} \tag{4.32}$$

It suffices to show that for every $t_1 > 0$ there exists a unique solution $\rho \in C^1([0,t_1]\times[0,1];(0,+\infty))$ of the initial-boundary value problem (4.27) with

$$\frac{\partial \rho}{\partial t}(t,x) + \lambda(W(t))\frac{\partial \rho}{\partial x}(t,x) = 0, \text{ for } (t,x) \in [0,t_1]\times[0,1] \tag{4.33}$$

$$W(t) = \int_0^1 u(t,x)dx, \text{ for } t \in [0,t_1] \tag{4.34}$$

$$\rho(t,0) = \rho_s \exp(b(t)), \text{ for } t \in [0,t_1] \tag{4.35}$$

which also satisfies the following estimate for all $t \in [0,t_1]$, $x \in [0,1]$:

$$\rho_{min} \leq \rho(t,x) \leq \rho_{max} \tag{4.36}$$

Let $L_{\rho_0} \leq \max\{|\rho_0'(x)|: x \in [0,1]\}$ be the Lipschitz constant for $\rho_0$ and $L_{\tilde{b}} \leq \max\left\{\left|\frac{d\tilde{b}}{dt}(t)\right|: t \in [0,t_1+1]\right\}$ be the Lipschitz constant for $\tilde{b}$ on $[0,t_1+1]$. We first show that for every $T \in \left(0, \left(1+L_\lambda\left(L_{\rho_0}+\frac{L_{\tilde{b}}}{v_{min}}\right)\right)^{-1}\right]$ there exists a unique solution $\rho \in C^1([0,T]\times[0,1];(0,+\infty))$ of the initial-boundary value problem (4.27), (4.29), (4.30), (4.31) with $t_1$ replaced by $T$, which satisfies estimate (4.36) for all $t \in [0,T]$, $x \in [0,1]$. Moreover, the Lipschitz constant $L_{\rho[t]}$ for the function $\rho[t]$ satisfies:

$$L_{\rho[t]} \leq \max\left(L_{\rho_0}, \frac{L_{\tilde{b}}}{v_{min}}\right), \text{ for all } t \in [0,T] \tag{4.37}$$

Indeed, if we show the above implications then we can construct step-by-step the solution of the initial-boundary value problem (4.27), (4.29), (4.30), (4.31), first on the interval $[0,T]$, then on the interval $[T,2T]$ and so on, with $T = \left(1+L_\lambda\left(\max\left(L_{\rho_0}, \frac{L_{\tilde{b}}}{v_{min}}\right)+\frac{L_{\tilde{b}}}{v_{min}}\right)\right)^{-1}$ and we cover the interval $[0,t_1]$.

Let arbitrary $T \in \left(0, \left(1+L_\lambda\left(L_{\rho_0}+\frac{L_{\tilde{b}}}{v_{min}}\right)\right)^{-1}\right]$ be given. Consider the operator $G: S \to S$ with

$$S := \left\{ v \in C^0([0,T]) : v_{\min} \leq \min_{t \in [0,T]}(v(t)) \leq \max_{t \in [0,T]}(v(t)) \leq v_{\max} \right\} \quad (4.38)$$

which maps the function $v \in S$ to the function $Gv = \bar{v} \in S$ defined by

$$\bar{v}(t) = \lambda \left( \int_0^1 \rho_v(t,x) dx \right), \text{ for } t \in [0,T] \quad (4.39)$$

where $\rho_v : [0,T] \times [0,1] \to \Re$ is the function defined by the equation

$$\rho_v(t,x) = \begin{cases} \rho_0 \left( x - \int_0^t v(s)ds \right) & \text{if } 1 \geq x > \int_0^t v(s)ds \\ \tilde{b}(t_0(t,x;v)) & \text{if } 0 \leq x \leq \int_0^t v(s)ds \end{cases} \quad (4.40)$$

where $t_0(t,x;v) \in [0,t]$ is the unique solution of the equation

$$x = \int_{t_0(t,x;v)}^t v(s)ds \text{ for all } (t,x) \in \Omega_v := \left\{ (t,x) \in [0,T] \times [0,1], x \leq \int_0^t v(s)ds \right\} \quad (4.41)$$

Notice that the mapping $\Omega_v \ni (t,x) \to t_0(t,x;v)$ is continuous with $|t_0(t,x;v) - t_0(t,y;v)| \leq \frac{1}{v_{\min}}|x-y|$ and $|t_0(t,x;v) - t_0(\tau,x;v)| \leq \frac{v_{\max}}{v_{\min}}|t-\tau|$ for all $(t,x) \in \Omega_v$, $(\tau,x) \in \Omega_v$, $(t,y) \in \Omega_v$. Moreover, when $x = \int_0^t v(s)ds \leq 1$ then $t_0(t,x;v) = 0$. Therefore, by virtue of the compatibility condition $\tilde{b}(0) = \rho_0(0)$, the mapping $\rho_v : [0,T] \times [0,1] \to \Re$ defined by (4.40) is continuous and satisfies the estimate $\rho_{\min} \leq \rho_v(t,x) \leq \rho_{\max}$ for all $t \in [0,T]$, $x \in [0,1]$. Moreover, due to the compatibility condition $\tilde{b}(0) = \rho_0(0)$, for every $t \in [0,T]$ the function $\rho_v[t]$ is Lipschitz on $[0,1]$ with Lipschitz constant $L_{\rho_v[t]} \leq \max\left( L_{\rho_0}, \frac{L_{\tilde{b}}}{v_{\min}} \right)$, where $L_{\rho_0} \leq \max\{|\rho_0'(x)| : x \in [0,1]\}$ is the Lipschitz constant for $\rho_0$ and $L_{\tilde{b}} \leq \max\left\{ \left|\frac{d\tilde{b}}{dt}(t)\right| : t \in [0, t_1+1] \right\}$ is the Lipschitz constant for $\tilde{b}$ on $[0, t_1+1]$.

Let two arbitrary functions $v, w \in S$ be given. Using (4.32), (4.39) we obtain

$$\|Gv - Gw\|_\infty \leq L_\lambda \max\{|\rho_v(t,x) - \rho_w(t,x)| : (t,x) \in [0,T] \times [0,1]\} \quad (4.42)$$

When $(t,x) \in \Omega_v$ and $(t,x) \in \Omega_w$, definition (4.41) implies that $\int_{t_0(t,x;w)}^t (w(s) - v(s))ds = \int_{t_0(t,x;v)}^{t_0(t,x;w)} v(s)ds$. Consequently, we obtain $T\|v-w\|_\infty \geq v_{\min}|t_0(t,x;w) - t_0(t,x;v)|$ which in conjunction with (4.40) gives

$$|\rho_v(t,x) - \rho_w(t,x)| = |\tilde{b}(t_0(t,x;v)) - \tilde{b}(t_0(t,x;w))| \leq L_{\tilde{b}} \frac{T}{v_{\min}} \|v-w\|_\infty \quad (4.43)$$

When $(t,x) \notin \Omega_v$ and $(t,x) \notin \Omega_w$, definition (4.40) implies that

$$\left|\rho_v(t,x) - \rho_w(t,x)\right| = \left|\rho_0\left(x - \int_0^t v(s)ds\right) - \rho_0\left(x - \int_0^t w(s)ds\right)\right|$$
$$\leq L_{\rho_0}\left|\int_0^t v(s)ds - \int_0^t v(s)ds\right| \leq L_{\rho_0} t\|v-w\|_\infty \leq L_{\rho_0} T\|v-w\|_\infty \qquad (4.44)$$

When $(t,x) \notin \Omega_v$ and $(t,x) \in \Omega_w$, i.e., when $\int_0^t w(s)ds \geq x > \int_0^t v(s)ds$ then we get

$$t\|v-w\|_\infty \geq x - \int_0^t v(s)ds > 0 \qquad (4.45)$$

Moreover, the inequality $\int_0^t w(s)ds \geq x > \int_0^t v(s)ds$ implies that $x - \int_{t_0(t,x;w)}^t v(s)ds > \int_0^t v(s)ds - \int_{t_0(t,x;w)}^t v(s)ds = \int_0^{t_0(t,x;w)} v(s)ds$ which combined with $x = \int_{t_0(t,x;w)}^t w(s)ds = \int_{t_0(t,x;w)}^t (w(s)-v(s))ds + \int_{t_0(t,x;w)}^t v(s)ds$ (a consequence of definition (4.41)) gives $\int_{t_0(t,x;w)}^t (w(s)-v(s))ds > \int_0^{t_0(t,x;w)} v(s)ds$. The previous inequality implies that

$$t\|v-w\|_\infty \geq v_{\min} t_0(t,x;w) \qquad (4.46)$$

Using (4.45), (4.46), the compatibility condition $\tilde{b}(0) = \rho_0(0)$ and definition (4.40), we obtain

$$\left|\rho_v(t,x) - \rho_w(t,x)\right| = \left|\rho_0\left(x - \int_0^t v(s)ds\right) - \tilde{b}(t_0(t,x;w))\right|$$
$$\leq \left|\rho_0\left(x - \int_0^t v(s)ds\right) - \rho_0(0)\right| + \left|\tilde{b}(0) - \tilde{b}(t_0(t,x;w))\right|$$
$$\leq L_{\rho_0} t\|v-w\|_\infty + L_{\tilde{b}} t_0(t,x;w) \qquad (4.47)$$
$$\leq L_{\rho_0} t\|v-w\|_\infty + L_{\tilde{b}} \frac{t}{v_{\min}}\|v-w\|_\infty$$
$$\leq T\left(L_{\rho_0} + \frac{L_{\tilde{b}}}{v_{\min}}\right)\|v-w\|_\infty$$

A similar estimate holds for the case $(t,x) \in \Omega_v$ and $(t,x) \notin \Omega_w$. Thus by combining (4.42), (4.43), (4.44) and (4.37) we get:

$$\|Gv - Gw\|_\infty \leq TL_\lambda \left(L_{\rho_0} + \frac{L_{\tilde{b}}}{v_{\min}}\right)\|v-w\|_\infty \qquad (4.48)$$

Since $T \in \left(0, \left(1 + L_\lambda\left(L_{\rho_0} + \frac{L_{\tilde{b}}}{v_{\min}}\right)\right)^{-1}\right]$, estimate (4.48) implies that the operator $G: S \to S$ is a contraction. Consequently, by virtue of Banach's fixed point theorem, there exists a unique $v \in S$ with $v = Gv$. It follows from (4.29), (4.39), (4.40), (4.41) and the compatibility conditions $\rho_s \exp(b(0)) = \rho_0(0)$, $\dot{b}(0) = -\lambda \left( \int_0^1 \rho_0(x)dx \right) \frac{\rho_0'(0)}{\rho_0(0)}$ that $\rho_v \in C^1([0,T]\times[0,1];(0,+\infty))$ is a solution of the initial-boundary value problem (4.27), (4.29), (4.30), (4.31) with $t_1$ replaced by $T$, which satisfies estimate (4.36) for all $t \in [0,T]$, $x \in [0,1]$.

Uniqueness follows from the fact that any solution of the initial-boundary value problem (4.27), (4.29), (4.30), (4.31) with $t_1$ replaced by $T$, necessarily gives a function $v \in S$ with $v = Gv$. The proof is complete. ◁

We are now ready to give the proof of Theorem 3.1.

**Proof of Theorem 3.1:** Every solution of the initial-boundary value problem (3.1), (3.2), (3.3), (3.4), (3.5) is a solution of the initial-boundary value problem (2.1), (2.2), (2.3) with $v(t,x) = \lambda(W(t))$ for $t \geq 0$, $x \in [0,1]$. Estimates (2.4), (2.5) in conjunction with estimate (4.28) imply estimates (3.6), (3.7). The proof is complete. ◁

## 5. Concluding Remarks

The present work provided a thorough study of stability of the 1-D continuity equation, which appears in many conservation laws. We have considered the velocity to be a distributed input and we have also considered boundary disturbances. Stability estimates are provided in all $L^p$ state norms with $p > 1$, including the case $p = +\infty$ (sup norm). However, in our Input-to-State Stability estimates, the gain and overshoot coefficients depend on the velocity. Moreover, the logarithmic norm of the state appears instead of the usual norm.

As remarked in the Introduction, the obtained results can be used in the stability analysis of larger models that contain the continuity equation. In the present paper, our results were used in a straightforward way for the stability analysis of non-local, nonlinear manufacturing models under feedback control. Working similarly, the obtained results can be used for the stability analysis of non-local traffic models.